\documentclass[12pt]{article}

\usepackage[T2A]{fontenc}
\usepackage[utf8]{inputenc}

\usepackage{amsfonts}
\usepackage{amssymb}
\usepackage{amsmath}
\usepackage{amsthm}

\def \le {\leqslant}
\def \ge {\geqslant}

\newtheorem{theorem}{Theorem}

\begin{document}

\begin{Large}
 \centerline{\bf Diophantine approximations with positive integers:}
\centerline{\bf some remarks}
\end{Large}
\vskip+0.5cm
\centerline{by {\bf Nikolay Moshchevitin}}
\vskip+0.5cm
\begin{small}
 {\bf Abstract.}\, We  
give some comments on our recent results related to W.M. Schmidt's conjecture and Diophantine exponents. 
\end{small}
\vskip+0.5cm
   
 This short communication is a supplement to our papers \cite{moshe,Mopositi}.

We consider a pair of  real numbers $\Theta =(\theta^1,\theta^2)$.
We are interested in small values of the linear form
$$
||\theta^1 m_1+\theta^2m_2
||$$
in positive integers $m_1, m_2$.
Put
$$
\psi (t) 
=\psi_{\Theta } (t) =\min_{m_1,m_2 \in \mathbb{Z} ,\,\, 0< \max (|m_1|,|m_2|) \le t}
||m_1\theta^1+m_2\theta^2||,
$$
$$
\psi^*(t) 
=\psi^*_{\Theta } (t) =\min_{x\in \mathbb{Z},\,\, 0< x \le t}
\max_{j=1,2}||x\theta^j||
$$
and
$$
\psi_+(t) 
=\psi_{+:\Theta } (t) =\min_{m_1,m_2 \in \mathbb{Z}_+,\,\, 0< \max (m_1,m_2) \le t}
||m_1\theta^1+m_2\theta^2||.
$$
Recall the definitions of Diophantine exponents
$$
\omega = \omega (\Theta) =
\sup\{ \gamma:\,\,\,
\liminf_{t \to \infty} t^\gamma \psi_{\Theta} (t) < \infty\},
$$ 
$$
\hat{\omega} = \hat{\omega} (\Theta) =
\sup\{ \gamma:\,\,\,
\limsup_{t \to \infty} t^\gamma \psi_{\Theta} (t) < \infty\}
$$
and
$$
\omega^* = \omega^* (\Theta) =
\sup\{ \gamma:\,\,\,
\liminf_{t \to \infty} t^\gamma \psi^*_{\Theta} (t) < \infty\},
$$ 
We introduce Diophantine exponents
$$
\omega_+ = \omega_+ (\Theta) =
\sup\{ \gamma:\,\,\,
\liminf_{t \to \infty} t^\gamma \psi_{+;\Theta} (t) < \infty\},
$$ 
and
$$
\hat{\omega}_+ = \hat{\omega}_{+} (\Theta) =
\sup\{ \gamma:\,\,\,
\limsup_{t \to \infty} t^\gamma \psi_{+;\Theta} (t) < \infty\}.
$$

\section{W.M. Schmidt's theorem and its extensions}

Put
$$
\phi = \frac{1+\sqrt{5}}{2} = 1.618^+.
$$
In  1976 
W.M. Schmidt \cite{SCH} proved the following theorem.
\vskip+0.5cm
\begin{theorem}[W.M. Schmidt]
Let real numbers $\theta^1_1,\theta^2$
be linearly independent over $\mathbb{Z}$
together with 1. Then there exists a sequence of integer two-dimensional vectors
 $(x_1(i), x_2(i))$
 such that

 1.\,\, $x_1(i), x_2(i) > 0$;

2.\,\, $||\theta^1x_1(i)+\theta^2x_2(i) ||\cdot (\max \{x_1(i),x_2(i)\})^\phi \to 0$ as $ i\to +\infty$.
\label{shpos}
\end{theorem}
\vskip+0.5cm
In fact W.M. Schmidt proved (see discussion in \cite{BUKR}) that  for $\theta^1,\theta^2$ under consideration
 one has the inequality
\begin{equation}\label{knipers}
 \omega_{+} \ge \max \left( \frac{\hat{\omega}}{\hat{\omega}-1}; \hat{\omega} - 1 +\frac{\hat{\omega}}{\omega}\right)
\end{equation}
from  which
we immediately deduce
$$
 \omega_{+} (\Theta)
\ge \phi
 .
$$
From Schmidt's argument one can easily see that for $\theta^1,\theta^2$ linearly independent together with 1 one has 
\begin{equation}\label{knipers}
 \hat{\omega}_{+} \ge   \frac{{\omega}}{{\omega}-1}.
\end{equation}
We would like to note here that Thurnheer (see Theorem 2 from \cite{aa1989}) showed that for $\theta^1,\theta^2$ linearly independent together with 1
in the case
\begin{equation}\label{knee}
\frac{1}{2}\le \omega^* = \omega^*(\Theta) \le 1 
\end{equation}
one has
\begin{equation}\label{knipers2}
 {\omega}_{+} \ge \frac{\omega^*+1}{4\omega^*} +\sqrt{\left(\frac{\omega^*+1}{4\omega^*} \right)^2+1}.
\end{equation}
(inequality  \ref{knipers2} is a particular case  of a general result obtained by Thurnheer).

A lower bound for $\omega_+$ in terms of $\omega  $ was obtained by the author in \cite{moshe}. 
It was based on the original Schmidt's 
argument from \cite{SCH}. However the choice of parameters in  \cite{moshe} was not optimal. Here we explain
the optimal choice.
From Schmidt's proof 
and Jarn\'{\i}k's result 
$$
\omega \ge \hat{\omega}(\hat{\omega}-1)
$$
(see \cite{J} and a recent paper \cite{L})
one can easily see that
\begin{equation}\label{zera}
\omega_{+} \ge
\max 
\left\{ g:\,\,\,
\max_{y,z \ge 1:\,\, y^{\hat{\omega}-1}\le z\le y^{\omega/\hat{\omega}}}
\,
\max_{y^{-\omega}\le x\le z^{-\hat{\omega}}}\,\,
\min\left(
x^{1-g}z^{-g}; x y^{-1}z^{g+1}
 \right) \le 1 
\right\}
.
\end{equation}
This inequality immediately follows from Schmidt's argument, see Lemma  1 and Lemma 2 from \cite{moshe}.
The right hand side of (\ref{zera}) can be easily calculated.
We divide the set
$$
\hbox{\got A} =
\left\{ (\omega, \hat{\omega}) \in \mathbb{R}^2:\,\,\, 
\hat{\omega}\ge 2,\,\, \omega \ge \hat{\omega}(\hat{\omega} -1) \right\}
$$
of all admissible values of $(\omega,\hat{\omega})$ into two parts:
$$
\hbox{\got A} =
\hbox{\got A}_1\cup \hbox{\got A}_2,
$$ 
$$
\hbox{\got A}_1 =
\left\{ (\omega, \hat{\omega}) \in \mathbb{R}^2:\,\,\,
2\le
\hat{\omega}\le \phi^2,\,\, \omega \ge \frac{\hat{\omega}(\hat{\omega}-1)}{3\hat{\omega}-\hat{\omega}^2 -1}\right\},
$$
$$
\hbox{\got A}_2 = \hbox{\got A} \setminus \hbox{\got A}_1
%\left\{ (\omega, \hat{\omega} \in \mathbb{R}^2:\,\,\,
%{\hat{\omega}(\hat{\omega}-1)}\le \omega \le 
%\frac{\hat{\omega}(\hat{\omega}-1)}{3\hat{\omega}-\hat{\omega}^2 -1}\right\},
.$$
If $(\omega,\hat{\omega}) \in \hbox{\got A}_1$
then
$$
\omega_+ \ge G(\omega) = \frac{1}{2}\left(\frac{\omega+1}{\omega}+\sqrt{\left(\frac{\omega+1}{\omega}\right)^2+4}\right)
$$
(the function $ G(\omega)$ on the right hand side decreases from $G(2) = 2$ to $G(+\infty ) = \phi$).
If 
$(\omega,\hat{\omega}) \in \hbox{\got A}_2$ then
\begin{equation}\label{ooa}
\omega_+ \ge \hat{\omega} - 1+\frac{\hat{\omega}}{\omega}
\end{equation}
So we get the following result.

\begin{theorem}
Let real numbers $\theta^1_1,\theta^2$
be linearly independent over $\mathbb{Z}$
together with 1. Then  
 $$
\omega_+ \ge
\max\left(\frac{1}{2}\left(\frac{\omega+1}{\omega}+\sqrt{\left(\frac{\omega+1}{\omega}\right)^2+4}\right);
\hat{\omega} - 1+\frac{\hat{\omega}}{\omega}\right)
.$$
\end{theorem}

This theorem  gives
the best bound in terms of $\omega, \hat{\omega}$ which one can deduce from Schmidt's argument from \cite{SCH}.

\section{About counterexample to W.M. Schmidt's conjecture
}

In the paper \cite{SCH} W.M. Schmidt wrote that he did not know if the exponent $\phi$ in Theorem \ref{shpos}
may be replaced by a lagrer constant. At that time he was not
able even to rule a possibility
that there exists an infinite sequence $(x_1(i),  x_2(i))\in \mathbb{Z}^2$ with condition 1. and such that
\begin{equation}\label{opip}
||\theta^1x_1(i)+\theta^2x_2(i) ||\cdot (\max \{x_1(i),x_2(i)\})^2 \le c(\Theta)
\end{equation}
with some large positive $c(\Theta)$.
Later in \cite{SCH1} he conjectured that the exponent $\phi$ may be replaced by any exponent of the form
$2 -\varepsilon, \varepsilon >0$ and wrote that probably such a result should be obtained
by analytical tools.
 It happened that this conjecture is not true. 
 In  \cite{Mopositi} the author proved the following result.
\begin{theorem}\label{moshepos1}
Let $\sigma = 1.94696^+$ be the largest real root of the equation
$
 x^4 - 2x^2-4x+1=0.
$ There exist real numbers
 $\theta^1,\theta^2$
such that they are linearly independent over $\mathbb{Z}$
together with 1 and for every  integer vector
$(m_1,m_2)\in\mathbb{Z}^2$ with $m_1 , m_2 \ge 0$  and
$\max (m_1,m_2) \ge 2^{200}$
one has
$$
||m_1\theta^1+m_2\theta^2||\ge \frac{1}{2^{300}(\max (m_1,m_2))^\sigma}.
$$
\end{theorem}

Here we should note that for the numbers constucted in Theorem \ref{moshepos1} one has
$$
\omega = \frac{(\sigma+1)^2(\sigma^2-1)}{4\sigma} = 3.1103^+ ,
\,\,\,
\hat{\omega} = \frac{(\sigma+1)^2}{2\sigma} = 2.2302^+ .
$$
So $(\omega, \hat{\omega}) \in \hbox{\got A}_2$ and the inequality (\ref{ooa}) gives
$$
\omega_+ \ge \frac{\sigma +2}{\sigma^2-1} = 1.413^+.
$$
However from the proof of Theorem \ref{moshepos1} (see \cite{Mopositi}) it is clear
that for the numbers constructed one has $\omega_+ = \sigma = 1.94696^+$.

\end{document}